\documentclass{jnmp}
\def\emline#1#2#3#4#5#6{
       \put(#1,#2){\special{em:moveto}}
      \put(#4,#5){\special{em:lineto}}}
\def\newpic#1{}

\makeatletter
\DeclareRobustCommand{\primfrac}[1]{%
  \PackageWarning{amsmath}{%
Foreign command \@backslashchar#1; %
\protect\frac\space or \protect\genfrac\space should be used instead%
  }
  \global\@xp\let\csname#1\@xp\endcsname\csname @@#1\endcsname
  \csname#1\endcsname
}
\makeatother
\begin{document}

\renewcommand{\evenhead}{I G Korepanov}
\renewcommand{\oddhead}{Invariants of PL Manifolds from Metrized Simplicial Complexes}

\thispagestyle{empty}


\FirstPageHead{8}{2}{2001}
{\pageref{Korepanov-firstpage}--\pageref{Korepanov-lastpage}}{Letter}

\copyrightnote{2001}{I G Korepanov}

\Name{Invariants of PL Manifolds from Metrized Simplicial Complexes.
Three-Dimensional Case}\label{Korepanov-firstpage}

\Author{Igor G KOREPANOV}

\Address{South Ural State University, 76 Lenin ave., Chelyabinsk 454080, Russia\\
E-mail: kig@susu.ac.ru}

\Date{Received October 6, 2000; Revised January 16, 2001; Accepted  January 17, 2001}

\begin{abstract}
\noindent
An invariant of three-dimen\-sional orientable manifolds is built on the
base of a solution of pentagon equation expressed in terms of metric
characteristics of Euclidean tetrahedra.
\end{abstract}

\section*{Introduction}

The motivation for work presented in this Letter is, from the side of
mathematical physics, in developing some ideas concerning building
of a topological field theory from a variant of ``Regge calculus''.
As is well-known, Regge~\cite{regge} proposed a discretization of
space-time in the form of its triangulation and assigning lengths
to the edges of such triangulation. In order to construct an analog of
functional integral in such theory, either a sum is taken over ``all''
triangulations or some sort of equivalence of different triangulations
is established.

The idea of such equivalence is well known in pure mathematics.
Mathematically, we want to put some numeric characteristics in
correspondence to a PL (piecewise-linear) manifold. The usual way of
doing this is, roughly speaking, as follows.
First, describe the manifold in algebraic or combinatorial terms~---
we have already done it as soon as we have chosen a triangulation.
Such a description can be done,
as a rule, in numerous different ways, but those different ways
can often be obtained from one another by using a sequence of some
simple re-buildings, or ``moves''. In the case of triangulation,
such moves affect only a few neighboring simplices from their total
maybe very big number.

Then, one could try to find an algebraic expression which could be
put in correspondence to some ``local'' part of the manifold, e.g.,
to a cluster of neighboring simplices, such that it would remain
similar to itself in some sense under the mentioned moves. Finally,
one could try to construct a ``global'' expression out of the ``local''
ones, in conformity with their algebraic structure, and find some
way to extract manifold invariants from such an expression.

A typical case of realization of the stated program is the building of
three-dimen\-sional manifold invariants out of quantum $6j$-symbols.
As key property of $6j$-symbols one can take the fact that they satisfy
the {\em pentagon equation\/} which is depicted in a natural way
as the equality of two diagrams, the first containing two tetrahedra with
a common base, while the second --- three tetrahedra occupying the
same domain in a Euclidean space, see Figure~\ref{fig 2-3} below.
Diagrams of such sort can be also introduced for the space of any
dimension~$n$ (the left- and right-hand sides must form together, at
least from the combinatorial viewpoint, the boundary of an
$(n+1)$-simplex). It seems however that any direct analogs of quantum
$6j$-symbols for higher-dimen\-sional manifolds are rather hard to find.

We would like to propose some other algebraic expressions that obey
a relation which, too, deserves the name of pentagon equation, because
the picture for it is the same. In constructing our expressions we assume
that the tetrahedra lie in a usual Euclidean space and thus possess
metric characteristics such as edge lengths, dihedral angles and
volumes. Our invariant is a certain expression made of those values.
Thus, it may be thought of as produced by some version of
Regge calculus. In order that our invariant be well defined, we will
assume that our PL manifold satisfies some additional requirements
including {\em orientability}.

The experience of the theory of discrete integrable models shows that
equations that are depicted by the same diagram turn out ultimately
to be closely connected, even if the diagram seems at first to have
completely different meanings. Besides, a connection of our expressions
with usual $6j$-symbols is suggested by Justin Roberts' work~\cite{JR}
where he explains how a metric tetrahedron appears in the quasiclassical
limit from $6j$-symbols corresponding to $SU(2)$ group (this was first
discovered by Ponzano and Regge~\cite{PR}, but not proved rigorously).
Here the quasiclassics is understood as tending of the irreducible
representations' dimensions to infinity. In this sense, our invariant looks
``quasiclassical'' too, but it is worth mentioning that,
again, the theory of integrable models teaches that the relations between
quantum and classical models are much richer than just classical models
being a limiting case of quantum ones, and in fact quantum models are
sometimes studied using ``classical'' considerations.

It looks plausible that our constructions can be generalized to
higher-dimen\-sional manifolds. Thus, the aim of the present Letter
is not only in introducing still more invariants of three-dimen\-sional
manifolds but in elaborating the necessary technical devices,
starting from this simplest case.

Below, in Section~\ref{sec local} we recall the derivation of the ``local''
formula from paper~\cite{dobavka}. This formula contains the partial
derivative of ``defect angle'' around an edge common for three
tetrahedra in the length  of that edge, taken in the neighborhood of the
flat case (when the whole cluster of tetrahedra can be imbedded into
a 3-dimen\-sional Euclidean space). Generalization (globalization)
of this formula onto the case of simplicial complexes having
{\em many\/} tetrahedra requires some technical work and occupies
Sections \ref{sec vzaim}--\ref{sec ves}.

In Section~\ref{sec vzaim} we introduce a matrix~$A$ of partial
derivatives of {\em all\/} defect angles (we call them simply
``curvatures'') in {\em all\/} edge lengths. Somewhat unexpectedly,
this matrix turns out to be {\em symmetric}: $A=A^{\rm T}$.
In Section~\ref{sec dobavl} we investigate matrix~$A$ from another
point of view: it is strongly degenerate, and we get first results in
globalizing the formula for $2\to 3$ (2 tetrahedra to 3 tetrahedra)
moves from Section~\ref{sec local}
by using matrix~$A$'s minor of the highest rank. It is also in
Section~\ref{sec dobavl} that we begin using the orientability
of the manifold. In Section~\ref{sec forma}
we continue our technical work and construct a differential form
not depending on the choice of the mentioned minor and behaving very nicely
under $2\leftrightarrow 3$ and $1\leftrightarrow 4$ moves.

\newpage

In Section~\ref{sec ves} we divide this differential form by some
``standard'' differential form and obtain a number that does not depend on
any edge length! Thus the explicit formula appears for the
invariant of a three-manifold.

To demonstrate the efficiency of our formula, we calculate
in Section~\ref{sec primery} our invariant in the two simplest
examples, namely, for the sphere~$S^3$ and the projective space~$RP^3$.
Finally, we discuss our results and their possible generalizations in
Section~\ref{sec discus}.

\section{The local formula}
\label{sec local}

In this section we recall the derivation of the formula from~\cite{dobavka}
that can be treated as a sort of pentagon equation involving five tetrahedra
in a three-dimen\-sional space. It corresponds in a natural way to replacing
a cluster of two Euclidean tetrahedra with the cluster of three tetrahedra
that covers the same 3-domain, or a $2\to 3$ {\em Pachner move},
as in Figure~\ref{fig 2-3}.

\begin{figure}[t]

\vspace*{5mm}
\begin{center}
\unitlength=1.00mm
\special{em:linewidth 0.4pt}
\linethickness{0.4pt}
\begin{picture}(121.00,47.00)
\emline{5.00}{30.00}{1}{20.00}{5.00}{2}
\emline{20.00}{5.00}{3}{35.00}{25.00}{4}
\emline{20.00}{5.00}{5}{50.00}{30.00}{6}
\emline{50.00}{30.00}{7}{35.00}{25.00}{8}
\emline{35.00}{25.00}{9}{5.00}{30.00}{10}
\put(4.00,30.00){\makebox(0,0)[rc]{$A$}}
\special{em:linewidth 0.2pt}
\linethickness{0.2pt}
\emline{5.00}{30.00}{11}{7.00}{30.00}{12}
\emline{9.00}{30.00}{13}{11.00}{30.00}{14}
\emline{13.00}{30.00}{15}{15.00}{30.00}{16}
\emline{17.00}{30.00}{17}{19.00}{30.00}{18}
\emline{21.00}{30.00}{19}{23.00}{30.00}{20}
\emline{25.00}{30.00}{21}{27.00}{30.00}{22}
\emline{29.00}{30.00}{23}{31.00}{30.00}{24}
\emline{33.00}{30.00}{25}{35.00}{30.00}{26}
\emline{37.00}{30.00}{27}{39.00}{30.00}{28}
\emline{41.00}{30.00}{29}{43.00}{30.00}{30}
\emline{45.00}{30.00}{31}{47.00}{30.00}{32}
\emline{49.00}{30.00}{97}{50.00}{30.00}{98}
\put(35.00,24.50){\makebox(0,0)[lt]{$B$}}
\put(51.00,30.00){\makebox(0,0)[lc]{$C$}}
\put(20.00,48.00){\makebox(0,0)[cb]{$D$}}
\put(20.00,4.00){\makebox(0,0)[ct]{$E$}}
\special{em:linewidth 0.4pt}
\linethickness{0.4pt}
\put(59.00,30.00){\vector(1,0){7.00}}
\emline{75.00}{30.00}{35}{90.00}{5.00}{36}
\emline{90.00}{5.00}{37}{105.00}{25.00}{38}
\emline{90.00}{5.00}{39}{120.00}{30.00}{40}
\emline{120.00}{30.00}{41}{105.00}{25.00}{42}
\emline{105.00}{25.00}{43}{75.00}{30.00}{44}
\put(74.00,30.00){\makebox(0,0)[rc]{$A$}}
\special{em:linewidth 0.2pt}
\linethickness{0.2pt}
\emline{75.00}{30.00}{45}{77.00}{30.00}{46}
\emline{79.00}{30.00}{47}{81.00}{30.00}{48}
\emline{83.00}{30.00}{49}{85.00}{30.00}{50}
\emline{87.00}{30.00}{51}{89.00}{30.00}{52}
\emline{91.00}{30.00}{53}{93.00}{30.00}{54}
\emline{95.00}{30.00}{55}{97.00}{30.00}{56}
\emline{99.00}{30.00}{57}{101.00}{30.00}{58}
\emline{103.00}{30.00}{59}{105.00}{30.00}{60}
\emline{107.00}{30.00}{61}{109.00}{30.00}{62}
\emline{111.00}{30.00}{63}{113.00}{30.00}{64}
\emline{115.00}{30.00}{65}{117.00}{30.00}{66}
\emline{119.00}{30.00}{67}{120.00}{30.00}{68}
\put(105.00,24.50){\makebox(0,0)[lt]{$B$}}
\put(121.00,30.00){\makebox(0,0)[lc]{$C$}}
\put(90.00,48.00){\makebox(0,0)[cb]{$D$}}
\put(90.00,4.00){\makebox(0,0)[ct]{$E$}}
\emline{90.00}{5.00}{69}{90.00}{7.00}{70}
\emline{90.00}{9.00}{71}{90.00}{11.00}{72}
\emline{90.00}{13.00}{73}{90.00}{15.00}{74}
\emline{90.00}{17.00}{75}{90.00}{19.00}{76}
\emline{90.00}{21.00}{77}{90.00}{23.00}{78}
\emline{90.00}{25.00}{79}{90.00}{27.00}{80}
\emline{90.00}{29.00}{81}{90.00}{31.00}{82}
\emline{90.00}{33.00}{83}{90.00}{35.00}{84}
\emline{90.00}{37.00}{85}{90.00}{39.00}{86}
\emline{90.00}{41.00}{87}{90.00}{43.00}{88}
\emline{90.00}{45.00}{89}{90.00}{47.00}{90}
\special{em:linewidth 0.4pt}
\linethickness{0.4pt}
\emline{75.00}{30.00}{93}{90.00}{47.00}{94}
\emline{90.00}{47.00}{95}{120.00}{30.00}{96}
\emline{105.00}{25.00}{97}{90.00}{47.00}{98}
\emline{20.00}{47.00}{99}{50.00}{30.00}{100}
\emline{35.00}{25.00}{101}{20.00}{47.00}{102}
\emline{20.00}{47.00}{103}{5.00}{30.00}{104}
\end{picture}
\end{center}
\vspace{-6mm}

\caption{A $2\to 3$ Pachner move}\label{fig 2-3}

\vspace{-2mm}

\end{figure}

Consider five points $A$, $B$, $C$, $D$ and $E$ in the three-dimen\-sional
Euclidean space. There exist ten distances between them, which we will
denote as $l_{AB}$, $l_{AC}$ and so on.

Let us fix all the distances except $l_{AB}$ and  $l_{DE}$. Then,
$l_{AB}$ and  $l_{DE}$ satisfy one constraint (Cayley--Menger equation)
which we can, using arguments like those in~\cite{3s-t}, represent
in the following differential form:
\begin{equation}
\left| l_{AB} \,dl_{AB} \over V_{\overline D}\, V_{\overline E} \right| =
\left| l_{DE} \,dl_{DE} \over V_{\overline A}\, V_{\overline B} \right|,
\label{eq AB-DE}
\end{equation}
where, say, $V_{\overline A}$ denotes the volume of tetrahedron $\overline A$, that is
one with vertices $B$, $C$, $D$ and $E$ (and {\em without\/}~$A$).

Let us consider the dihedral angles at the edge $DE$ --- the common edge
for tetrahedra $\overline A$, $\overline B$ and $\overline C$. Namely, denote
\[
\angle BDEC \stackrel{\rm def}{=} \alpha, \qquad
\angle CDEA \stackrel{\rm def}{=} \beta, \qquad
\angle ADEB \stackrel{\rm def}{=} \gamma.
\]
We have:
\begin{equation}
0=d(\alpha+\beta+\gamma) = \frac{\partial \gamma}{\partial l_{AB}} \,dl_{AB} +
\frac{\partial (\alpha+\beta+\gamma)}{\partial l_{DE}} \,dl_{DE}.
\end{equation}

According to~\cite[formula~(11)]{3s-t},
\begin{equation}
\left|\frac{\partial \gamma}{\partial l_{AB}}\right|={1\over 6}\left| l_{AB}\,l_{DE}\over
V_{\overline C}\right|.
\label{eq 3}
\end{equation}
Denote also
\begin{equation}
\alpha+\beta+\gamma \stackrel{\rm def}{=} 2\pi-\omega_{DE},
\label{eq dgamma}
\end{equation}
where $\omega_{DE}$ is the ``defect angle'' around edge~$DE$. The formulas
(\ref{eq AB-DE})--(\ref{eq dgamma}) together yield
\begin{equation}
\left|{1\over V_{\overline D}\, V_{\overline E}}\,  \right| =
{1\over 6}\, \left| {l_{DE}^2\over V_{\overline A}\, V_{\overline B}\, V_{\overline C}}
\left( \frac{\partial\omega_{DE}}{\partial l_{DE}} \right)^{-1} \right|.
\label{eq local pa}
\end{equation}

\medskip

\noindent
{\bf Remark.} This can be also written by means of the following integral in the
length of the edge~$DE$ ``redundant'' for the tetrahedra $\overline D$ and
$\overline E$:
\begin{equation}
\left|{1\over V_{\overline D}\, V_{\overline E}} \right| =
{1\over 6}\, \left| \int {\delta(\omega_{DE})\, l_{DE}^2 \,dl_{DE}
\over V_{\overline A}\, V_{\overline B}\, V_{\overline C} } \right| ,
\label{eq delta}
\end{equation}
with the integral taken over a neighborhood of the value of~$l_{DE}$
corresponding to the flat space ($\omega_{DE}=0$); $\delta$ is the Dirac
delta function. However, the straightforward attempt to globalize
formula~(\ref{eq delta}) runs into diverging integrals, and the right
``global'' formulae (see, e.g.,~(\ref{eq **6})) will have volumes raised
into the power $(-1/2)$ rather than~$(-1)$.

\section{Reciprocity theorems for lengths and defect angles}
\label{sec vzaim}

In this Section we will do some of the technical work mentioned in the
Introduction. Consider a finite simplicial complex made of tetrahedra
and their faces (of dimensions 2, 1 and~0). If the contrary is not stated
explicitely, we assume that every 2-face belongs to boundaries of exactly
two tetrahedra lying at its different sides (thus, the corresponding
PL-manifold as a whole has no boundary).

Assign to each edge of the complex a {\em length\/}, say length $l_a$
to the edge~$a$. Consider a cluster of all tetrahedra containing the
edge~$a$. The edge lengths in this cluster may happen to be consistent in
such way that the whole cluster can be put into a Euclidean 3-space.
Generally, however, there is an obstacle called the {\em defect angle
$\omega_a$} corresponding to edge~$a$ which we define {\em up to a
multiple of $2\pi$} by the equality
\[
\omega_a \equiv -\alpha-\beta-\cdots-\eta \pmod{2\pi},
\]
where $\alpha,\beta,\ldots,\eta$ are the proper dihedral angles of
the tetrahedra.

Now we will consider the partial derivatives like $\partial\omega_a
/ \partial l_b$ which are taken with the fixed lengths of all edges
except~$b$. It must be clear from the preceding paragraph that such 
a~partial derivative may be nonzero only if the edges $a$ and~$b$ belong
to a single tetrahedron.

\medskip

\noindent
{\bf Theorem 1 (local reciprocity theorem).}
{\it Let a tetrahedron in the Euclidean space be given, $a$ and $b$ being
its two edges (they can lie on skew, intersecting or coinciding
straight lines),
$l_a$ and $l_b$ being their lengths, and $\varphi_a$ and $\varphi_b$ ---
dihedral angles at those edges. Then
\begin{equation}
\frac{\partial\varphi_a}{\partial l_b}=\frac{\partial \varphi_b}{\partial l_a}.
\label{eq local}
\end{equation}
}

\medskip

\noindent
{\bf Proof.} The case of coinciding edges $a$ and $b$ is trivial.

The case of skew edges: both l.h.s.\ and r.h.s.\ of (\ref{eq local})
equal $(1/6)l_a l_b/V$, where $V$ is the volume of tetrahedron
(compare formula~(\ref{eq 3})).

\begin{figure}[t]
\vspace*{5mm}

\begin{center}
\unitlength=1.00mm
\special{em:linewidth 0.4pt}
\linethickness{0.4pt}
\begin{picture}(51.00,51.00)
\emline{5.00}{30.00}{1}{20.00}{5.00}{2}
\emline{20.00}{5.00}{3}{35.00}{25.00}{4}
\emline{35.00}{25.00}{5}{20.00}{50.00}{6}
\emline{20.00}{50.00}{7}{5.00}{30.00}{8}
\emline{20.00}{5.00}{9}{50.00}{30.00}{10}
\emline{50.00}{30.00}{11}{35.00}{25.00}{12}
\emline{50.00}{30.00}{13}{20.00}{50.00}{14}
\emline{35.00}{25.00}{15}{5.00}{30.00}{16}
\put(4.00,30.00){\makebox(0,0)[rc]{$A$}}
\special{em:linewidth 0.2pt}
\linethickness{0.2pt}
\emline{5.00}{30.00}{17}{7.00}{30.00}{18}
\emline{9.00}{30.00}{19}{11.00}{30.00}{20}
\emline{13.00}{30.00}{21}{15.00}{30.00}{22}
\emline{17.00}{30.00}{23}{19.00}{30.00}{24}
\emline{21.00}{30.00}{25}{23.00}{30.00}{26}
\emline{25.00}{30.00}{27}{27.00}{30.00}{28}
\emline{29.00}{30.00}{29}{31.00}{30.00}{30}
\emline{33.00}{30.00}{31}{35.00}{30.00}{32}
\emline{37.00}{30.00}{33}{39.00}{30.00}{34}
\emline{41.00}{30.00}{35}{43.00}{30.00}{36}
\emline{45.00}{30.00}{37}{47.00}{30.00}{38}
\emline{49.00}{30.00}{39}{50.00}{30.00}{40}
\put(35.00,24.50){\makebox(0,0)[lt]{$B$}}
\put(51.00,30.00){\makebox(0,0)[lc]{$C$}}
\put(20.00,51.00){\makebox(0,0)[cb]{$D$}}
\put(20.00,4.00){\makebox(0,0)[ct]{$D'$}}
\emline{20.00}{5.00}{41}{20.00}{7.00}{42}
\emline{20.00}{9.00}{43}{20.00}{11.00}{44}
\emline{20.00}{13.00}{45}{20.00}{15.00}{46}
\emline{20.00}{17.00}{47}{20.00}{19.00}{48}
\emline{20.00}{21.00}{49}{20.00}{23.00}{50}
\emline{20.00}{25.00}{51}{20.00}{27.00}{52}
\emline{20.00}{29.00}{53}{20.00}{31.00}{54}
\emline{20.00}{33.00}{55}{20.00}{35.00}{56}
\emline{20.00}{37.00}{57}{20.00}{39.00}{58}
\emline{20.00}{41.00}{59}{20.00}{43.00}{60}
\emline{20.00}{45.00}{61}{20.00}{47.00}{62}
\emline{20.00}{49.00}{63}{20.00}{50.00}{64}
\put(25.00,27.50){\makebox(0,0)[cb]{$a$}}
\put(40.00,27.50){\makebox(0,0)[cb]{$b$}}
\end{picture}
\end{center}
\vspace{-6mm}

\caption{To the proof of Theorem~1}
\label{fig ABCDD'}

\vspace{-7mm}

\end{figure}
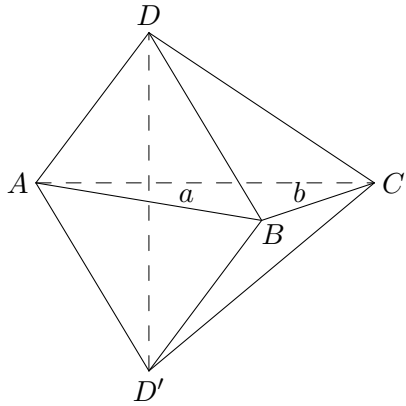

The case of intersecting edges. Let $a$ and~$b$ be, respectively, edges
$AB$ and $BC$ in a~tetrahedron $ABCD$. Consider also the mirror image
$ABCD'$ of tetrahedron $ABCD$ with respect to the plane $ABC$, see
Figure~\ref{fig ABCDD'}.

Now we can calculate, say, $\partial\varphi_b/\partial l_a$ in the
following way. Assuming that all the edge lengths except $a$ and $DD'$
in Figure~\ref{fig ABCDD'} are fixed, let us calculate first
$\partial l_{DD'}/\partial l_a$. We will get, similarly to
formula~(\ref{eq AB-DE}) (and using the same notations like $V_{\overline A}$):
\begin{equation}
\frac{\partial l_{DD'}}{\partial l_a}=-\frac{l_a}{l_{DD'}}
\frac{V_{\overline A} V_{\overline B}}{V_{ABCD}^2}.
\label{eq ***1}
\end{equation}

Next, from tetrahedron $BCDD'$ (in other words --- tetrahedron~$\overline A$)
we can get (compare formula~(\ref{eq 3})):
\begin{equation}
2 \frac{\partial \varphi_b}{\partial l_{DD'}}=
\frac 16 \frac{l_b l_{DD'}}{V_{\overline A}}.
\label{eq ***2}
\end{equation}

It follows from (\ref{eq ***1}) and (\ref{eq ***2}) that
\[
\frac{\partial \varphi_b}{\partial l_a}=-\frac{1}{12}
\frac{l_a l_b V_{\overline B}}{V_{ABCD}^2}.
\]

Clearly, the result will be the same for $\partial\varphi_a/\partial l_b$.
The theorem is proved.

\medskip

\noindent
{\bf Theorem 2 (global reciprocity theorem).}
{\it Let a complex be given of the type described in the beginning of this
Section. Select in it two edges $a$ and $b$. Then
\begin{equation}
\left.\frac{\partial \omega_a}{\partial l_b}\right|_{\hbox
{\scriptsize
\rm $l_c$ are constant for $c\ne b$ } }
=
\left.\frac{\partial \omega_b}{\partial l_a}\right|_{\mbox
{\scriptsize
\rm $l_c$ are constant for $c\ne a$ } }.
\label{eq rec gen}
\end{equation}
}

\medskip

\noindent
{\bf Proof.} The equality (\ref{eq rec gen})
follows from the fact that the l.h.s.\ of (\ref{eq rec gen}) is the sum of
values of the type $-\partial\varphi_a^{(k)}/\partial l_b$, where $k$
numbers the tetrahedra containing {\em both edges\/} $a$ and~$b$, and
$\varphi_a^{(k)}$ is the dihedral angle in such tetrahedron at edge~$a$.
As for the r.h.s.\ of (\ref{eq rec gen}), it is the sum of similar terms
but with interchanged $a\leftrightarrow b$, and these sums are equal due to
the local reciprocity theorem. The theorem is proved.

\medskip

\noindent
{\bf Addition to Theorem~2.}
The equality~(\ref{eq rec gen}) remains valid if we change the definitions
of defect angles in the following way: select any subset in the set of
tetrahedra of the complex, and assume {\em all\/} dihedral angles in those
tetrahedra to be {\em negative}.

\medskip

\noindent
{\bf Proof} follows immediately from an obvious modification of the
proof of Theorem~2.

\medskip

Introduce the matrix
\begin{equation}
A=\left( \frac{\partial \omega_j}{\partial l_k} \right),
\label{matr A}
\end{equation}
where $j$ and $k$ run through all the edges of the complex.
Matrix~$A$ is thus {\em symmetric\/}: $A=A^{\rm T}$.

\section{A quantity good for $\pbf{2\to 3}$ moves}
\label{sec dobavl}

Throughout the rest of this Letter, we will be considering ``metrized''
simplicial 3-comp\-le\-xes (with lengths assigned to their edges) of the type
described in the beginning of
Section~\ref{sec vzaim} with the following two additional
constraints: {\em the corresponding PL manifold must be orientable,
and the lengths are such that the corresponding polyhedron
can be put into the 3-dimen\-sional Euclidean space}~$R^3$.
This will be understood as follows: we identify all {\em vertices\/} of the
complex with points in~$R^3$. Thus, all edges acquire Euclidean
lengths, and every tetrahedron gets embedded in~$R^3$. It is important
to note that we do permit any self-intersections of the obtained
Euclidean tetrahedra.

We say that such lengths form a {\em permitted length configuration\/}.
At the same time,
we will consider {\em any infinitesimal deformations\/} of lengths which
can thus draw the complex out of the Euclidean space or, in other words,
produce some infinitesimal defect angles around edges. For brevity, we
will sometimes call those defect angles ``curvatures'' (as an exception from
this rule, we will soon be considering a situation where {\em one\/} ``new''
edge can take any value, but if we remove it the complex fits again
into Euclidean space).

The accurate definition of infinitesimal defect angles shows why we require
the orientability of the manifold. Fix a consistent orientation of all
tetrahedra in the complex. When we map the complex into a Euclidean space
(which we suppose to have its own fixed orientation), some of the tetrahedra
preserve their orientation while the others change it. We will define
the defect angle around a given edge as the algebraic sum of (interior)
dihedral angles in adjoining tetrahedra taken with the sign~$-$ for
the tetrahedra that do not change their orientation and with the
sign~$+$ for the rest of them. Such definition ensures that the defect
angles in a complex mapped into Euclidean space will be zero (in absense
of infinitesimal deformations), and we will be using the Addition to
Theorem~2 exactly in such situation.

The matrix~$A$ given by (\ref{matr A}) is usually strongly degenerate,
see examples below in Section~\ref{sec primery}. It follows from
the fact that $A$ is symmetric and standard theorems in linear algebra
that there exists a {\em diagonal\/} nondegenerate submatrix $A|_{\cal C}$
of~$A$ of sizes $\mbox{rank}\; A \times \mbox{rank}\; A$. This means that we can
choose a subset~$\cal C$ in the set of all edges, and leave only those
rows and columns in~$A$ that correspond (both rows and columns) to
edges from~$\cal C$.

Now we will make the considerations of the preceding paragraph more
precise in the following way. Matrix~$A$ depends on a chosen permitted
length configuration. We will be dealing
with the ranks of $A$ and its submatrices for a {\em generic\/}
permitted configuration.
Accordingly, below we denote by $\cal C$ a chosen
subset of $L$ edges for which $A|_{\cal C}$ is nondegenerate in the
general position, where $L$ equals $\mbox{rank}\; A$ again in the general position.

The rest of edges form the subset that we will denote~$\overline{\cal C}$.

\medskip

\noindent
{\bf Lemma 1.} {\it  The form $\bigwedge\limits_{i\in\overline{\cal C}} dl_i$, i.e.\ the exterior
product of differentials of all edge lengths from $\overline{\cal C}$, is 
nondegenerate in a generic point of the algebraic variety consisting of
all permitted length configurations.}

\medskip

\noindent
{\bf Proof.} The lemma can be reformulated as follows: for any set of
length differentials $dl_i$ of edges in $\overline{\cal C}$ one can find
such length differentials of edges in $\cal C$ that {\em all\/}
infinitesimal curvatures $d\omega$ will equal zero. Now, it follows
immediately from the nondegeneracy of matrix $A|_{\cal C}=(\partial
\omega_j / \partial l_k)|_{\cal C}$ that one can always find such
differentials of lengths in $\cal C$ that all the infinitesimal curvatures
{\em around edges in $\cal C$} will be zero. But any other
infinitesimal curvature {\em is linearly dependent\/} upon the curvatures
in $\cal C$ and thus vanishes as well. The lemma is proved.

\medskip

Consider two adjacent tetrahedra (with a common 2-face) in the complex
and perform the operation of replacing them with three tetrahedra, as in
Section~\ref{sec local}. In doing so, we add a new edge ($DE$
in Figure~\ref{fig 2-3}) of length~$l_{\rm new}$. Denote $\tilde l_{\rm new}
\stackrel{\rm def}{=} l_{\rm new}-l_{\rm new}^{(0)}$, where
$l_{\rm new}^{(0)}$ is
such value of $l_{\rm new}$ where the curvature $\omega_{\rm new}$ around
the new edge is exactly zero. Then
\begin{equation}
d\tilde l_{\rm new} = dl_{\rm new} -a_1\, dl_1 -\cdots -a_N\, dl_N,
\label{eq *0}
\end{equation}
where $N$ is the number of edges before adding the new one, and
\[
a_k=\frac{\partial l_{\rm new}^{(0)}}{\partial l_k}=
\left.\frac{\partial l_{\rm new}}{\partial l_k}\right|_{\hbox{\scriptsize
\begin{tabular}{l}
$\omega_{\rm new}=0$\\
$dl_1=\cdots=dl_N=0$
\end{tabular}}}
= -\frac{\partial \omega_{\rm new} / \partial l_k}
{\partial \omega_{\rm new} / \partial l_{\rm new} }.
\]
It is clear that $d\omega_{\rm new}$ depends only on $d\tilde l_{\rm new}$
and does not depend on $dl_1,\ldots,\,dl_N$:
\begin{equation}
d\omega_{\rm new} = \frac{\partial \omega_{\rm new}}{\partial l_{\rm new}}\,
d\tilde l_{\rm new},
\label{eq *1}
\end{equation}
where we have taken into account the obvious equality
\begin{equation}
\left.\frac{\partial}{\partial \tilde l_{\rm new}}\right|_{dl_1=\cdots=dl_N=0} =
\left.\frac{\partial}{\partial l_{\rm new}}\right|_{dl_1=\cdots=dl_N=0}.
\label{eq *2}
\end{equation}
Then, when $d\tilde l_{\rm new}=0$, the definition of matrix~$A$ works:
\begin{equation}
\left(\begin{array}{c}
d\omega_1\\ \vdots\\ d\omega_N\end{array}\right) =
A\left(\begin{array}{c} dl_1\\ \vdots\\ dl_N\end{array}\right).
\label{eq *3}
\end{equation}
Combining (\ref{eq *1}), (\ref{eq *2}) and (\ref{eq *3}) we can write:
\begin{equation}
\left(\begin{array}{c} d\omega_{\rm new}\\ d\omega_1\\ \vdots\cr d\omega_N
\end{array}\right) =
\left(\begin{array}{cc} \partial \omega_{\rm new} / \partial l_{\rm new} &
0 \  \cdots \ 0\vspace{1mm}\\
\begin{array}{c}
\partial \omega_1 / \partial l_{\rm new}\\
 \vdots \\
  \partial \omega_N / \partial l_{\rm new}\end{array}
 &
 \hbox{\Huge{$A$}}\end{array}\right)
\left(\begin{array}{c}
 d\tilde l_{\rm new}\\ dl_1\\ \vdots\\ dl_N\end{array}\right).
\label{eq *4}
\end{equation}

We are willing to construct the new matrix $A_{\rm new}$ of sizes
$(N+1)\times (N+1)$ that will link, like matrix~$A$ did, the differentials
of lengths and curvatures, but for the complex with the added edge.
Combining (\ref{eq *0}) and~(\ref{eq *4}) we get:
\begin{equation}
A_{\rm new} =
\left(\begin{array}{cc}
\partial \omega_{\rm new} / \partial l_{\rm new} &
 0 \ \cdots \ 0 \vspace{1mm}\\
\begin{array}{c} \partial \omega_1 / \partial l_{\rm new}\\
 \vdots \\   \partial \omega_N / \partial l_{\rm new} \end{array}
& \hbox{\Huge{$A$}}\end{array}\right)
\left(
\begin{array}{cccc}
1 & -a_1 & \cdots & -a_N \\
           &  1   &        &      \\
           &      & \ddots &  \hbox{\Huge 0}\\
           & \hbox{\Huge 0}  &        &   1 \end{array}\right).
\label{eq *5}
\end{equation}

Let ${\cal C}_{\rm new}$ be the subset of the set of edges obtained
by adding the ``new'' edge to $\cal C$.
It can be seen from formula (\ref{eq *5}) that we can take
$A_{\rm new}|_{{\cal C}_{\rm new}}$ for a submatrix of $A_{\rm new}$
having the same rank as~$A_{\rm new}$. To be exact, it follows from
(\ref{eq *5}) that
\begin{equation}
\det \left(A_{\rm new}|_{{\cal C}_{\rm new}}\right) =
\frac{\partial \omega_{\rm new}}{\partial l_{\rm new}} \det \left(A|_{\cal C}\right).
\label{eq *6}
\end{equation}
Note that $\partial \omega_{\rm new} / \partial l_{\rm new}$ has been
calculated in Section~\ref{sec local}, see formula (\ref{eq local pa}).
With this taken into account, formula~(\ref{eq *6}) shows that the
following theorem is valid.

\medskip

\noindent
{\bf Theorem 3.} {\it The expression
\begin{equation}
\frac{f}{\prod\limits_{\hbox{\scriptsize \rm over all edges}} l^2}
\prod_{\hbox{\scriptsize \rm over all tetrahedra}} 6V,
\label{eq *7}
\end{equation}
where
\begin{equation}
f \stackrel{\rm def}{=} \det \left(A|_{\cal C}\right),
\label{eq *8}
\end{equation}
and $l$'s and $V$'s are of course lengths and volumes,
does not change or changes only its sign under performing a
$2\to 3$ Pachner move in such way that a new edge is added to
the subset $\cal C$.}

\medskip

In order to construct out of (\ref{eq *7}) an invariant of a PL manifold,
we still have to get rid of the dependence of our construction on the
concrete choice of subset~$\cal C$ and also make our formulas describe
not only $2\leftrightarrow 3$ moves (adding or removing an edge, as in
Section~\ref{sec local}), but also $1\leftrightarrow 4$, when a new
{\em vertex\/} is added to the complex or removed. This is what we will
do next.

\section{Differential forms and Pachner moves}
\label{sec forma}

In this Section we are going to consider the following column vectors
of differentials:
\begin{itemize}
\item
$d{\pbf{l}}$ --- the column of differentials of lengths for edges from
$\cal C$, i.e.
\[
d{\pbf{l}}=\left(\begin{array}{c} \vdots\\ dl_i\\ \vdots\end{array}\right),\quad
i\in {\cal C};
\]
\item
$d{\pbf{k}}$ --- the column of differentials of lengths for edges from
$\overline{\cal C}$, i.e.
\[
d{\pbf{k}}=\left(\begin{array}{l} \vdots\\ dl_i\\ \vdots\end{array}\right),\quad
i\in \overline{\cal C};
\]
\item
$d{\pbf{\omega}}$ --- the column of differentials of curvatures around
edges from $\cal C$;
\item
$d{\pbf{\psi}}$ --- the column of differentials of curvatures around
edges from $\overline{\cal C}$.
\end{itemize}

\noindent
{\bf Lemma 2.} {\it Matrix $A$ (introduced in (\ref{matr A})), if written in a block form
corresponding to the above mentioned partitions of sets of differentials:
\begin{equation}
\left(\begin{array}{c} d{\pbf{\omega}}\\
 d{\pbf{\psi}}\end{array}\right)=
A\left(\begin{array}{c} d{\pbf{l}}\\
d{\pbf{k}}\end{array}\right),
\label{eq **0}
\end{equation}
has the following block structure:
\begin{equation}
A=\left(\begin{array}{cc} A|_{\cal C} & -(A|_{\cal C}){\pbf{a}}\vspace{1mm}\\
 -{\pbf{a}}^{\rm T}(A|_{\cal C}) & {\pbf{a}}^{\rm T}(A|_{\cal C}){\pbf{a}}\end{array}\right)
=\left(\begin{array}{cc} {\bf 1} & {\bf 0}\\ -{\pbf{a}}^{\rm T} & {\bf 0}\end{array}\right)
\left(\begin{array}{cc} A|_{\cal C} & {\bf 0} \\ {\bf 0} & {\bf 0}\end{array}\right)
\left(\begin{array}{cc} {\bf 1} & -{\pbf{a}} \\ {\bf 0} & {\bf 0}\end{array}\right),
\label{eq **.5}
\end{equation}
where $\pbf{a}$ is the matrix connecting $d{\pbf{k}}$ and 
$d{\pbf{l}}$
in the flat case:
\begin{equation}
d{\pbf{l}}|_{d{\pbf{\scriptstyle\omega}}=0} = {\pbf{a}}\, d{\pbf{k}}
\label{eq **.6}
\end{equation}
(recall that it follows from $d{\pbf{\omega}}=0$ that $d{\pbf{\psi}}=0$
as well); the superscript $\rm T$ means matrix transposing.}

\medskip

\noindent
{\bf Proof.} Introduce (in analogy with the situation of adding a ``new''
edge in Section~\ref{sec dobavl}) the following column of differentials:
\begin{equation}
d\tilde{\pbf{l}}=d{\pbf{l}}-{\pbf{a}}\, d{\pbf{k}}.
\label{eq **1}
\end{equation}
Then if $d\tilde{\pbf{l}}$ is zero, all curvatures vanish as well.
This can be written as
\begin{equation}
\left(\begin{array}{c} d{\pbf{\omega}}\\ d{\pbf{\psi}}\end{array}\right)=
\left(\begin{array}{cc} B & {\bf 0}\\ C & {\bf 0}\end{array}\right)
\left(\begin{array}{c} d\tilde{\pbf{l}}\\ d{\pbf{k}}\end{array}\right),
\label{eq **2}
\end{equation}
where $B$ and $C$ are some matrices, and it will be clear soon that
$B=A|_{\cal C}$.

On the other hand, (\ref{eq **1}) can be rewritten as
\begin{equation}
\left(\begin{array}{c} d\tilde{\pbf{l}}\\ d{\pbf{k}}\end{array}\right)=
\left(\begin{array}{cc} {\bf 1} & -{\pbf{a}}\\ {\bf 0} & {\bf 1}\end{array}\right)
\left(\begin{array}{c} d{\pbf{l}}\\ d{\pbf{k}}\end{array}\right).
\label{eq **3}
\end{equation}

Comparing (\ref{eq **2}) and (\ref{eq **3}) on the one hand, and
(\ref{eq **0}) on the other, we see that
\[
A=\left( \begin{array}{cc} B & -B{\pbf{a}}\\ C & -C{\pbf{a}}\end{array}\right).
\]

It is clear now that $B=A|_{\cal C}$, while the blocks in the
second row of matrix~$A$
are determined from the fact that it is symmetric
(Theorem~2 and the Addition to it). The lemma is proved.

\medskip

Now we are ready to investigate how $f=\det \left(A|_{\cal C}\right)$
changes under replacing $\cal C$ with a similar set ${\cal C}'$. Such a
${\cal C}'$ can be obtained like this: choose subsets ${\cal A}\in{\cal C}$
and ${\cal B}\in\overline{\cal C}$ with the same number of elements, and move
$\cal A$ from $\cal C$ to $\overline{\cal C}$, while $\cal B$ --- vice versa.
$\cal C$ is thus transformed into a set which we can take for~${\cal C}'$
($\cal A$ and $\cal B$ must be chosen, however, in such a way that
$\det \left(A|_{{\cal C}'}\right)$ be not identically zero).

Denote also
\[
f'=\det \left(A|_{{\cal C}'}\right).
\]

\medskip

\noindent
{\bf Lemma 3.} {\it 
\begin{equation}
\frac{f'}{f}=\left( \det \left(\vphantom{|}_{\cal A}|{\pbf{a}}|_{\cal B} \right) \right)^2,
\label{eq **3.5}
\end{equation}
where $\vphantom{|}_{\cal A}|{\pbf{a}}|_{\cal B}$
 means the square submatrix of~$\pbf{a}$ for which the rows
corresponding to edges from~$\cal A$ and the columns corresponding to
edges from~$\cal B$ are taken.}

\medskip

\noindent
{\bf Proof.} One can see from the form (\ref{eq **.5}) of matrix~$A$
that $A|_{{\cal C}'}$ can be expressed through~$A|_{\cal C}$. Namely,
\begin{equation}
A|_{{\cal C}'} = F^{\rm T} A|_{\cal C} F,
\label{eq **4}
\end{equation}
where the square matrix $F$ has the following block structure that
corresponds to set partition ${\cal C}={\cal C}\backslash{\cal A}\cup
{\cal A}$ for the rows and ${\cal C}'={\cal C}'\backslash{\cal B}\cup
{\cal B}$ for the columns:
\begin{equation}
F=\left(\begin{array}{cc} {\bf 1} & {\bf *}\\ {\bf 0} & -\vphantom{|}_{\cal A}|{\pbf{a}}|_{\cal B}
\end{array}\right)
\label{eq **5}
\end{equation}
(note that the {\em whole\/} matrix~$F$ is thus of the same size
as the upper left blocks e.g.~in formula~(\ref{eq **.5})).
Here the asterisk means some submatrix that we are not interested in.
By applying (\ref{eq **4}) and (\ref{eq **5}) the lemma is proved.

\medskip

It follows from the definition (\ref{eq **.6}) of matrix~$\pbf{a}$
connecting length differentials for zero curvatures that the determinant
in the r.h.s.\ of (\ref{eq **3.5}) is nothing but
\[
\det\left(\vphantom{|}_{\cal A}|{\pbf{a}}|_{\cal B}\right)=
 \pm \frac{\bigwedge\limits_{i\in{\cal A}} dl_i}
{\bigwedge\limits_{i\in{\cal B}} dl_i}=
\pm \frac{\bigwedge\limits_{i\in\overline{{\cal C}'}} dl_i}
{\bigwedge\limits_{i\in\overline{\cal C}} dl_i},
\]
where all $dl_i$ are taken as well for zero curvatures, i.e.\ they are
forms on the space of {\em permitted length configurations}. They are
nonzero due to Lemma~1. Hence, an important
conclusion follows: {\em the form
\[
\frac{1}{\sqrt{|f|}} \bigwedge\limits_{i\in\overline{\cal C}} dl_i,
\quad \mbox{or simply} \quad \frac{1}{\sqrt{|f|}}
\bigwedge\limits_{\overline{\cal C}} dl,
\]
does not change, to within its sign, with a different choice of~$\cal C$}.
This conclusion can be united with Theorem~3 in the following way.

\medskip

\noindent
{\bf Theorem 4.} {\it The following differential form on the variety of permitted length
configurations on the edges of complex (the degree of the form coincides
with the dimension of variety) remains unchanged under
$2\leftrightarrow 3$ Pachner moves and under a different choice of
subset $\cal C$ in the set of edges:
\begin{equation}
\left|
\frac{\prod\limits_{\hbox{\scriptsize \rm over all edges}} l
\; \bigwedge\limits_{\overline{\cal C}} dl}
{\sqrt{|f|\; \prod\limits_{\hbox{\scriptsize \rm over all tetrahedra}} 6V}}
\right|.
\label{eq **6}
\end{equation}}

Consider now a $1\to 4$ Pachner move. This means that a tetrahedron
$ABCD$ is replaced with four tetrahedra $ABCE$, $ABDE$, $ACDE$ and $BCDE$,
where $E$ is a new vertex added to the complex. We will assume that
edge~$DE$ is added to the set~$\cal C$, while edges $AE$, $BE$ and $CE$
are added to the set~$\overline{\cal C}$. To trace the changes in the
form~(\ref{eq **6}), we have to note that $f$ is multiplied by
$\partial \omega_{DE} / \partial l_{DE}$, and this partial derivative can
be calculated again from formula~(\ref{eq local pa}) (although we are now
in a somewhat different situation). As a result, (\ref{eq **6}) is
multiplied by
\begin{equation}
\left|\frac{l_{AE}\, dl_{AE} \wedge l_{BE}\, dl_{BE} \wedge
l_{CE}\, dl_{CE}}{6V_{ABCE}}
\right|.
\label{eq **7}
\end{equation}

Assume that all edge lengths are temporarily fixed except $l_{AE}$, $l_{BE}$
and $l_{CE}$, and the complex is put into the
3-dimen\-sional Euclidean space with a fixed coordinate system $Oxyz$
in such way that coordinates of all its vertices except~$E$ are fixed.
Then a simple trigonometry shows that (\ref{eq **7}) turns into
\begin{equation}
\left|
dx_{E} \wedge dy_{E} \wedge dz_{E}
\right|,
\label{eq **8}
\end{equation}
where, of course, $x_{E}$, $y_{E}$ and $z_{E}$ are Euclidean coordinates
of point~$E$.

\section{The invariant $\pbf{I}$}
\label{sec ves}

Let us select three vertices among the vertices of the complex and denote
them $A$, $B$ and~$C$. Draw the axes $x$, $y$ and $z$ of a Euclidean
system of coordinates in such way that~$A$ be the origin of coordinates,
$B$ lie on the $x$ axis and $C$ --- in the plane~$xAy$. When we vary the
lengths of edges of the complex in a ``permitted'' way, the coordinates
of vertices, namely $x_B$, $x_C$, $y_C$, $x_D$, $y_D$, $z_D, \ldots$
change, too.

\medskip

\noindent
{\bf Lemma 4.} {\it The form
\begin{equation}
\left| x_B^2\, dx_B \wedge dx_C \wedge y_C \, dy_C \wedge
\bigwedge_{\hbox{\scriptsize over remaining vertices}} dx \wedge
dy \wedge dz \right|
\label{eq 3 *}
\end{equation}
does not depend on the choice of $A$, $B$ and~$C$.}

\medskip

\noindent
{\bf Proof.} This simple fact can be proved in different ways. We prefer
to link it to the ideas of paper~\cite{3s-t}.

Let us begin with the case where there is only one vertex~$D$ in the
complex besides $A$, $B$ and~$C$. The already mentioned trigonometry
(see (\ref{eq **7}) and~(\ref{eq **8})) shows that with
fixed $A$, $B$ and~$C$
\[
|dx_D \wedge dy_D \wedge dz_D| = 
\frac{|l_{AD}\, dl_{AD} \wedge
l_{BD}\, dl_{BD} \wedge l_{CD}\, dl_{CD}|}{6 V_{ABCD}}.
\]
Then, it is also easy to show that
\begin{equation}
|x_B^2\, dx_B \wedge dx_C \wedge y_C \, dy_C| =
|l_{AB}\, dl_{AB} \wedge l_{AC}\, dl_{AC} \wedge l_{BC}\, dl_{BC}|.
\label{eq 3 *.}
\end{equation}
Thus, the whole expression~(\ref{eq 3 *}) is in this case equal to
\begin{equation}
\frac{|l_{AB}\, dl_{AB} \wedge \cdots \wedge l_{CD}\, dl_{CD}|}
{6V_{ABCD}},
\label{eq 3 **}
\end{equation}
with the exterior product in the enumerator taken over all edges of
tetrahedron~$ABCD$. It is clear that~(\ref{eq 3 **}) does not change
under all permutations of the set $\{A,B,C,D\}$.

Let now there be {\em five\/} vertices $A$, $B$, $C$, $D$ and~$E$
in the complex. Then it is easy to show that (\ref{eq 3 *}) turns into
\begin{equation}
\frac{|l_{AB}\, dl_{AB} \wedge \cdots \wedge l_{CE}\, dl_{CE}|}
{6V_{ABCD}\cdot 6V_{ABCE}},
\label{eq 3 ***}
\end{equation}
where the exterior product in the enumerator is taken over all edges
entering in {\em at least one\/} of tetrahedra $ACBD$ and~$ABCE$ (in
other words, the edge $DE$ is absent from~(\ref{eq 3 ***})). To conclude
the proof for five vertices it is enough to show that (\ref{eq 3 ***})
does not change under the permutation $C\leftrightarrow D$ (note that,
again, $A$, $B$ and $C$ can be interchanged freely because
of~(\ref{eq 3 *.})).

Under $C\leftrightarrow D$, the factor $V_{ABCE}$ in the denominator of
(\ref{eq 3 ***}) is replaced with $V_{ABDE}$, and $l_{CE}\, dl_{CE}$
in the enumerator is replaced with $l_{DE}\, dl_{DE}$.
Now it remains to apply the formula
\[
\left| \frac{l_{CE}\, dl_{CE}}{V_{ABCE}} \right| =
\left| \frac{l_{DE}\, dl_{DE}}{V_{ABDE}} \right|,
\]
compare~\cite[formula~(12)]{3s-t}.

Finally, if there are more than five vertices, the proof of the lemma
is obtained by obvious generalization of the above arguments. The lemma
is proved.

\medskip

Now let us note that the degree of the form~(\ref{eq 3 *}) increases
or decreases under moves $1\leftrightarrow 4$ in the same way as
the degree of the form~(\ref{eq **6}), and they both do not change under
$2\leftrightarrow 3$ moves. Thus, the {\em difference\/} of those
degrees is already a manifold invariant! Below, however, we will concentrate
on another invariant which is defined only for those manifolds that
satisfy the following assumption.

\medskip

\noindent
{\bf Assumption.} Let our PL manifold be such that the degrees of forms
(\ref{eq 3 *}) and (\ref{eq **6}) coincide and, moreover, they are
proportional in each point of the ``permitted'' variety.

\medskip

This Assumption seems to be satisfied for manifolds with finite fundamental
groups, see examples in Section~\ref{sec primery}.

Adopting this Assumption, we divide (\ref{eq **6}) by (\ref{eq 3 *}).
In general, we expect to get some function of the edge lengths in the
complex. It turns out, somewhat surprisingly, that it actually {\em does not
depend on those lengths\/} and is thus a constant depending on the manifold
only!

To see this, let us return to our arguments about the $1\to 4$ move from
the end of Section~\ref{sec forma}. When we add the point~$E$, the form
(\ref{eq **6}) is multiplied by the form (\ref{eq **7}) or, which is the
same, (\ref{eq **8}). Thus, the ratio $\hbox{(\ref{eq **6})}/
\hbox{(\ref{eq 3 *})}$, firstly, does not change and, secondly, is obviously
independent of the position of point~$E$ with respect to other vertices
of the complex, that is of lengths $l_{AE}$, $l_{BE}$ and $l_{CE}$.

The point $E$ has however equal rights in this respect with the other
vertices of the complex (any vertex can be eliminated, after some
preparatory moves, by a $4\to 1$ move. Thus, it can be regarded as added
to some complex by the reverse $1\to 4$ move). So, the ratio
$\hbox{(\ref{eq **6})}/\hbox{(\ref{eq 3 *})}$ (if it is well defined)
does not depend at all on how we place the vertices in the Euclidean space.

\medskip

\noindent
{\bf Definition.} {\it We define the invariant $I(M)$ for a given closed
oriented PL manifold $M$ by the formula
\[
I(M)=\frac{\mbox{\rm (\ref{eq **6})}}{\mbox{\rm (\ref{eq 3 *})}},
\]
provided the right-hand side is well defined.}

\section{Examples}
\label{sec primery}

In order to calculate the invariant~$I$ for the sphere~$S^3$,
it is enough to use its decomposition in 2 tetrahedra. In such way we get
a {\em pre-simplicial complex\/} (see e.g.~\cite{hilton wiley}) rather
than a simplicial complex, but our formulae remain valid for
this case as well.

There will be 4 vertices, say $A,B,C$
and $D$; 6~edges which {\em all\/} enter the set~$\overline{\cal C}$;
2~identical volumes~$V_{ABCD}$; and the value~$f$, due to the emptiness
of the set~$\cal C$, will be simply~$1$ (so, the rank of matrix
$A=(\partial\omega_j / \partial l_k)$ will be zero). Thus, the
formula~(\ref{eq **6}) will yield in this case exactly the
expression~(\ref{eq 3 **}). This means that
\[
I\left(S^3\right)=1.
\]

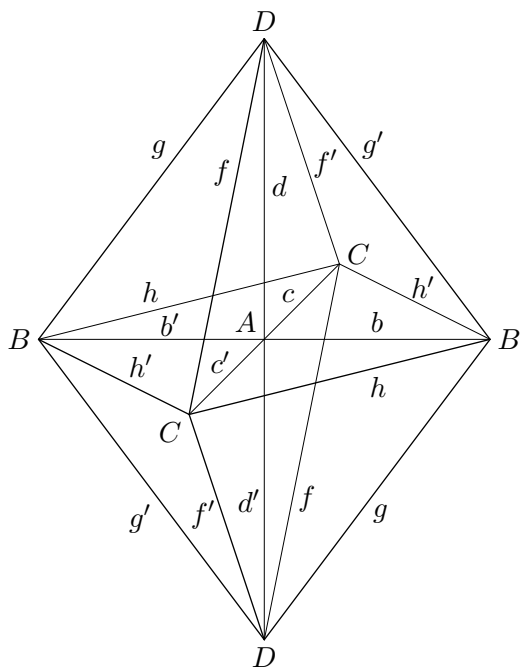
\begin{figure}[t]
\vspace*{5mm}

\begin{center}
\unitlength=1mm
\special{em:linewidth 0.5pt}
\linethickness{0.5pt}
\begin{picture}(66.00,86.00)
\emline{5.00}{45.00}{1}{25.00}{35.00}{2}
\emline{25.00}{35.00}{3}{65.00}{45.00}{4}
\emline{65.00}{45.00}{5}{35.00}{5.00}{6}
\emline{35.00}{5.00}{7}{25.00}{35.00}{8}
\emline{35.00}{5.00}{9}{5.00}{45.00}{10}
\emline{5.00}{45.00}{11}{35.00}{85.00}{12}
\emline{35.00}{85.00}{13}{25.00}{35.00}{14}
\emline{65.00}{45.00}{15}{35.00}{85.00}{16}
\special{em:linewidth 0.15pt}
\linethickness{0.15pt}
\emline{5.00}{45.00}{17}{45.00}{55.00}{18}
\emline{45.00}{55.00}{19}{65.00}{45.00}{20}
\emline{65.00}{45.00}{21}{5.00}{45.00}{22}
\emline{25.00}{35.00}{23}{45.00}{55.00}{24}
\emline{45.00}{55.00}{25}{35.00}{85.00}{26}
\emline{35.00}{85.00}{27}{35.00}{5.00}{28}
\emline{35.00}{5.00}{29}{45.00}{55.00}{30}
\put(22.00,69.00){\makebox(0,0)[rb]{$g$}}
\put(48.00,69.00){\makebox(0,0)[lb]{$g'$}}
\put(30.50,65.50){\makebox(0,0)[rb]{$f$}}
\put(41.50,66.50){\makebox(0,0)[lb]{$f'$}}
\put(36.00,64.00){\makebox(0,0)[lb]{$d$}}
\put(20.00,23.00){\makebox(0,0)[rt]{$g'$}}
\put(28.50,23.50){\makebox(0,0)[rt]{$f'$}}
\put(34.50,25.00){\makebox(0,0)[rt]{$d'$}}
\put(39.50,25.50){\makebox(0,0)[lt]{$f$}}
\put(49.50,23.00){\makebox(0,0)[lt]{$g$}}
\put(17.00,40.00){\makebox(0,0)[lb]{$h'$}}
\put(49.00,40.00){\makebox(0,0)[lt]{$h$}}
\put(54.50,50.50){\makebox(0,0)[lb]{$h'$}}
\put(21.00,50.00){\makebox(0,0)[rb]{$h$}}
\put(22.50,45.50){\makebox(0,0)[cb]{$b'$}}
\put(50.00,46.00){\makebox(0,0)[cb]{$b$}}
\put(39.00,50.00){\makebox(0,0)[rb]{$c$}}
\put(30.50,40.50){\makebox(0,0)[rb]{$c'$}}
\put(34.00,46.00){\makebox(0,0)[rb]{$A$}}
\put(4.00,45.00){\makebox(0,0)[rc]{$B$}}
\put(66.00,45.00){\makebox(0,0)[lc]{$B$}}
\put(46.00,55.00){\makebox(0,0)[lb]{$C$}}
\put(24.00,34.00){\makebox(0,0)[rt]{$C$}}
\put(35.00,4.00){\makebox(0,0)[ct]{$D$}}
\put(35.00,86.00){\makebox(0,0)[cb]{$D$}}
\end{picture}
\end{center}
\vspace{-6mm}

\caption{Triangulation for $RP^3$}
\label{fig RP3}
\end{figure}

Consider now the projective space~$RP^3$. For it, we take a triangulation
that has, again, 4 vertices $A$, $B$, $C$ and $D$ but now 12~edges
and 8~tetrahedra, see Figure~\ref{fig RP3}.
The edges will be distributed between sets $\cal C$ and $\overline{\cal C}$ as
follows: ${\cal C}=\{b,c,d,f,g,h\}$;
$\overline{\cal C}=\{b',c',d',f',g',h'\}$.

\medskip

\noindent
{\bf Remark.} Figure~\ref{fig RP3} represents the ``abstract triangulation''
of $RP^3$ and does {\em not\/} depict the imbedding of the vertices in
Euclidean space. Of course, such an imbedding cannot send, say,
vertex~$B$ in two different points.

\medskip

Most partial derivatives of curvatures in edge lengths in
Figure~\ref{fig RP3} are zero. We will explain this on the example of
$\partial\omega_d / \partial l_f$. There are exactly two tetrahedra
containing both edges $d$ and~$f$. The dihedral angles at edge~$d$
in those tetrahedra coincide in absolute value but have opposite signs,
and this remains so even when we vary~$l_f$ (in this case, it means that
when taking the derivative $\partial\omega_d / \partial l_f$ we fix
all lengths except $l_f$ so that $l_b=l_{b'}$, $l_c=l_{c'}$, $l_d=l_{d'}$,
$l_g=l_{g'}$ and $l_h=l_{h'}$ but allow $l_f$ to vary in a neighborhood
of $l_{f'}$).

Arguments of such sort show that the block $A|_{\cal C}$ of matrix
$A=(\partial\omega_i / \partial l_j)$ contains only 6 nonzero entries,
namely
\begin{equation}
\frac{\partial \omega_b}{\partial l_f},\ \
 \frac{\partial \omega_f}{\partial l_b},\ \
 \frac{\partial \omega_c}{\partial l_g},\ \
\frac{\partial \omega_g}{\partial l_c},\ \
 \frac{\partial \omega_d}{\partial l_h}\  \
\mbox{and} \ \
\frac{\partial \omega_h}{\partial l_d},
\label{eq der}
\end{equation}
so that $\det A|_{\cal C}$ coincides in absolute value with the
product of them all.

In order to calculate, say, $\partial\omega_b / \partial l_f$ we note
that there exist exactly two tetrahedra containing both edges
$b$ and~$f$, and one can see from Figure~\ref{fig RP3} that in these
tetrahedra the derivatives of dihedral angles at edge~$b$ in $l_f$
have the {\em same\/} sign. They are calculated by formulae of type
(\ref{eq 3}), so that we get
\begin{equation}
\left| \frac{\partial \omega_b}{\partial l_f} \right| = 2\cdot \left| 
\frac{l_b\, l_f}{6\, V_{ABCD}} \right|.
\label{eq mnozhitel 2}
\end{equation}

The rest of derivatives (\ref{eq der}) are calculated in a similar way,
and due to the multiplication by 2 in (\ref{eq mnozhitel 2}) and other
such formulas the form (\ref{eq **6}) turns out to be $1/8$ of the
similar form for $S^3$, which means that
\[
I\left(RP^3\right)=\frac{1}{8}.
\]

\section{Discussion}
\label{sec discus}

In this Letter we have only got some first results showing that
it is possible to construct, on the base of such quantities as edge
lengths, dihedral angles and volumes, an invariant which
can be calculated with no complications at least for the sphere~$S^3$
and projective space~$RP^3$. Actually, some calculations have been done
also for lens spaces $L(p,q)$, and they suggest that
\begin{equation}
I\bigl(L(p,q)\bigr)=\frac{1}{p^3}.
\label{eq lens}
\end{equation}
Thus, this version of invariant seems, somewhat regretfully, not to
depend on~$q$.

It must not be forgotten however that the presented version of invariant
is only the simplest one (and actually another invariant has been already
mentioned in Section~\ref{sec ves}, before the Assumption). An interesting
possibility is to map in~$R^3$ not a simplicial complex itself but its
{\em universal covering}, in case of a nontrivial fundamental group~$\pi_1$.

This idea will most likely be combined with the form (\ref{eq 3 *})
being replaced with some other ``standard'' differential form, of the
same degree as~(\ref{eq **6}). So, we will be led probably to richer
invariant structures than just a number like~(\ref{eq lens}).

The work~\cite{3s-t} suggests that invariants of the same type as in this
Letter can be constructed for higher-dimen\-sional manifolds as well.
This may be combined with other ways of generalization such as
the use of noncommutative ``lengths''.

Finally, the fact that the tetrahedron volumes enter in
formula~(\ref{eq **6}) raised in the power $(-1/2)$ shows that our
formulae are akin to the quasiclassical formulae suggested by Ponzano and
Regge~\cite{PR} and proved in a recent work by Justin Roberts~\cite{JR}
(where very interesting mathematical facts are presented related to
$6j$-symbols and Euclidean geometry).

\subsection*{Acknowledgements}
I am grateful to many people for scientific discussions, explaining
me topological and other subtleties and simply for encouraging letters.
Among them are: Rinat Kashaev, Greg Kuperberg, Marco Mackaay,
Sergei Matveev, Justin Roberts, Satoru Saito, and Jim Stasheff.
I also acknowledge the partial financial support from Russian Foundation
for Basic Research under grant no.~01-01-00059.


\strut\hfill

\label{Korepanov-lastpage}

\end{document}